\documentclass[reqno,twoside,12pt]{amsart}

\usepackage{amsfonts}

\newcommand{\spn}{\operatorname{span}}

\newcommand{\la}{\lambda}

\newcommand{\e}{\varepsilon}
\newcommand{\f}{\varphi}

\newcommand{\C}{\mathbb{C}}
\newcommand{\D}{\mathbb{D}}
\newcommand{\T}{\mathbb{T}}

\newcommand{\R}{\mathbb{R}}
\newcommand{\Z}{\mathbb{Z}}

\newcommand{\ci}[1]{_{ {}_{\scriptstyle #1}}}

\newtheorem{thm}{\hspace{\parindent}Theorem}[section]
\newtheorem{lm}[thm]{\hspace{\parindent}Lemma}
\newtheorem{cor}[thm]{\hspace{\parindent}Corollary}

\newcommand{\dist}{{\operatorname{dist}}}
\newcommand{\trace}{{\operatorname{trace}}}

\newcommand{\bE}{\mathbb{E}}

\begin{document}

\title
[Completely regular   stationary processes]
{Completely regular multivariate  stationary processes and the 
Muckenhoupt condition}
\author{S. Treil and A. Volberg}
\address{Department of Mathematics, Michigan State University,
East Lansing, Michigan, 48824}
\email[Treil]{treil@math.msu.edu}
\email[Volberg]{volberg@math.msu.edu}
\curraddr[Volberg]{Mathematical Sciences Research Institute,
1000 Centennial Drive, Berkeley, CA 94707-5070}

\thanks{Partially supported by the NSF
grant DMS 9622936, binational Israeli-USA grant BSF 00030, and
research program at MSRI in the Fall of 1997.}
\subjclass{42B20, 42A50, 47B35}

\begin{abstract}
We are going to give necessary and sufficient conditions for a multivariate stationary stochastic process to be completely regular. We also give the answer to a question of V.V. Peller concerning the spectral measure characterization of such processes.
\end{abstract}

\maketitle

\section{Introduction}
\label{intro}

In this paper we shall give a necessary and sufficient condition for a 
multivariate stationary stochastic process to be {\em completely 
regular}.  For the scalar case the description of completely regular processes 
was obtained by Helson an Sarason, see \cite{HSar,Sar}. Almost none of the scalar methods is available in the vector situation. The explanation is simple. Our problem will be reduced to verifying $L^2$ weighted inequalities for a certain integral operator. The weight will be a  matrix weight arising from the spectral measure of the process.  All the pointwise estimates  of integral operators become too crude for the vector valued case. For example, if a positive kernel is majorized by another one, and this second kernel gives the bounded operator in $L^2(\mu)$, then the original kernel obviously corresponds to a bounded operator in $L^2(\mu)$ too. But this is not the case if $\mu$ is a matrix measure even for scalar kernels. 

The study of prediction theory for multivariate stationary stochastic processes
was started by Kolmogorov and Wiener in the 50's, see, for example \cite{WM1}, \cite{WM2}, and \cite{MW}. It was later continued in works of I. Ibragimov, Yu. Rozanov, V. Solev, A. Yaglom, V. Peller, S. Khruschev, N.J. Young. An extensive  bibliography can be found in \cite{PKh} (for scalar processes) and in \cite{Pe} (for vector ones).
153a158,169

Let us recall that a multivariate stationary stochastic process with 
discrete time
is a sequence of $d$-tuples 
$x(n) = (x_{1}(n), x_{2}(n), \ldots , x_{d}(n))$, $n\in \Z$ of scalar 
random variables such that $\bE |x_{j}(n)|^{2} <\infty$ and the {\em 
correlation matrix} $Q(n,k)$
$$
Q(n,k)= \left\{ Q(n,k)_{i,j}\right\}_{1\le i,j\le d}:=
\left\{ \bE x_{i}(n) \overline{x_{j}(k)} \right\}_{1\le i,j\le d}
$$
depends only on the difference $n-k$; here $\bE$ denotes mathematical 
expectation.

It is well known (see \cite{Roz}) that there exists a matrix-valued non-negative 
measure $M$ on the unit circle $\T$ whose Fourier coefficients 
coincide with entries of the correlation matrix
$$
Q(n,k)= \widehat{M}(n-k)\,,\qquad n,k\in \Z\,.
$$
The measure $M$ is called the \emph{spectral measure} of the process $\{ 
x(n)\}_{n\in \Z}$.

The random variables $x_{j}(n)$ can be treated as elements of Hilbert 
space $L^{2}(\Omega, dP)$, where $\Omega$ is the probability space 
and $P$ is the probability, so $x(n)$ can be treated as elements of
the $\R^d$-valued $L^2$ space $L^{2}_{\R^d}(\Omega, dP)$
For a moment $n$ of time we can consider the past $\mathcal{X}_{n}$ 
and the future $\mathcal{X}^{n}$ of the process, which are defined as 
the subspaces
\begin{align*}
	\mathcal{X}_{n} & = \spn\left\{ x_{j}(k)\,:\, 1\le j \le d,\, 
	k<n\right\}
	\\
	\mathcal{X}^{n} & = \spn\left\{ x_{j}(k)\,:\, 1\le j \le d,\, 
	k\ge n\right\}
\end{align*}
of $L^{2}(\Omega, dP)$. 

A process is called {\em regular} if $\cap_{n\ge0}\mathcal X^{n} = 
\{0\}$. In this case (see \cite{Roz}) the spectral measure $M$ of the 
process is absolutely continuous with respect to Lebesgue measure. 
Let $W$ be the density of $M$ with respect to Lebesgue measure. The 
matrix-valued function $W$ is called the spectral density of the 
process. 

\medskip

A process $\{x(n)\}_{n\in\Z}$ is called {\em completely regular} if 
its past is asymptotically orthogonal to the future, namely if 
$$
\sup\left\{ | \bE (\xi {\eta}) |\,:\, \xi 
\in\mathcal{X}_{0},\, \eta\in \mathcal{X}^{n},\, \bE |\xi|^{2}\leq 
1,\, \bE |\eta|^{2}\leq 1 \right\} \longrightarrow 0\qquad \mathrm{as }\  n\to  
\infty\,.
$$
Of course, complete regularity implies regularity. 
If the process is Gaussian (i.e. all random variables $x_{j}(k)$ have 
normal distribution) then the complete regularity means simply that 
past and future are almost independent. The problem we are dealing 
with is to characterize completely regular processes in terms of 
spectral measure.

It is well known (see again \cite{Roz}) that if the process is 
completely regular, then its spectral measure is absolutely 
continuous, $dM =W dm$ where $dm$ is the normalized ($m(\T)=1$) Lebesgue 
measure on the unit circle $\T$. 

The reader is referred to  \cite{Roz} once more to see that there exists 
$d_{0}\leq d$ (the rank of the process) such that the spectral 
density $W(t)$ has rank $d_{0}$ for almost all $t\in\T$. If $d_{0} = 
d $ then the process $\{x(n)\}$ is said to be a {\em full rank}. 

The study of processes of arbitrary rank can be easily reduced to the 
study of the processes of full rank, see \cite{Ib}. So in this paper 
we shall consider only processes of full rank.

For the scalar case the description of completely regular processes 
was obtained by Helson an Sarason, see \cite{HSar,Sar}. To state their 
result we need a couple of definitions. 

Let us recall that a function $f$ on the unit circle $\T$ belongs to 
the space BMO (bounded mean oscillation) if
$$
\sup_{I} \frac{1}{|I|} \int_{I} |f-f_{I}| dm = \| 
f\|_{\mathrm{BMO}} <\infty\,;
$$
here $f_{I}$ denotes the mean value of $f$ on the interval $I$: 
$f_{I}:= |I|^{-1} \int_{I} f dm$ and the {\em supremum} is taken over 
all subarcs $I$ of $\T$. 

The space VMO (vanishing mean oscillation) consists of all function $f \in \mathrm{BMO}$ such that 
$$
\sup_{I} \frac{1}{|I|} \int_{I} |f-f_{I}| dm \longrightarrow 0\qquad 
\mathrm{as}\ |I|\to 0\, . 
$$

\begin{thm}[Helson, Sarason]
\label{He-Sar}
	Let $w$ be the spectral density of a scalar stationary 
	process. Then the process is completely regular if and only if $w$ 
	admits a representation
	$$
	w=|p|^{2}e^{\varphi}\,,
	$$
	where $p$ is a polynomial with roots on the unit circle $\T$ and 
	$\varphi$ is a real-valued function in VMO. 
\end{thm}

It was conjectured by V. Peller in \cite{Pe} that the same result 
holds for multivariate stationary processes. Namely he conjectured 
that a multivariate stationary process is completely regular if and 
only if its spectral density  $W$ admits the following representation
$$
W = P^{*} e^{\Phi} P, 
$$  
where $P$ is a polynomial matrix whose determinant has roots on $\T$ 
and the matrix function $\Phi = \Phi^{*}$ belongs VMO.

In this direction  he was able to prove the following theorem

\begin{thm}
	\label{t-Pe1}
	A multivariate stationary process is completely regular if and only 
	if its spectral density $W$ admits the  factorization
	$$
	W = P^{*} W_{1} P, 
	$$
	where $P$ is a polynomial matrix whose determinant has roots on $\T$ 
	and $W_{1}$ is the density of a completely regular stationary process 
	such that $W^{-1}_{1}\in L^{1}$. 
\end{thm}

\subsection{The main result}
\label{main-result}

Let us recall that a measure $\mu$ on the unit disk $\D$ is called 
Carleson if 
$$
\sup_{I} \mu(Q(I)) \le C\cdot |I|
$$
and is  called the vanishing Carleson measure 
if
$$
\limsup_{|I|\to 0}\mu(Q(I))/|I| = 0
$$
where {\em limsup} is taken over all subarcs $I$ of $T$. Here $Q(I)$ denotes the ``Carleson square'' for the arc $I$, 
$$
Q(I) = \{z\in \D\,:\, z/|z|\in I, 1-|I|\le |z| <1\}
$$ 

For a function $F$ on the unit circle let $F(\la)$, $\la\in \D$, denote 
its harmonic extension at the point $\la$. 

The main result of the paper is the following theorem.

\begin{thm}
	\label{t-main}
	Let  the density $W$  of a stationary process satisfy $W^{-1}\in 
	L^{1}$. 
	Then the the following are equivalent
	\begin{enumerate}
		\item  The process is completely regular; 
		\item $W^{-1}$ is the spectral density of a completely regular 
		process; 
		\item  ${\displaystyle \limsup_{|I|\to 0}\left\| 
		\Bigl(\frac{1}{|I|} \int_{I} W dm \Bigr)^{1/2} \Bigl(\frac{1}{|I|} 
		\int_{I} W^{-1} dm \Bigr)^{1/2} \right\| = 1 }$; here {\em 
		supremum} is taken over all subarcs $I$ of $\T$;
	
		\item  ${\displaystyle \limsup_{|\la|\to 1}\left\| 
		\Bigl(W(\la) \Bigr)^{1/2} \Bigl(W^{-1}(\la) \Bigr)^{1/2} \right\| = 
		1 }$, where $W(\la)$ and $W^{-1}(\la)$ are harmonic extensions of 
		functions $W\Bigm|  \T$ and $W^{-1}\Bigm |\T$ respectively at point 
		$\la\in\D$. 
	
	
		\item  ${\displaystyle \limsup_{|\la|\to 1}\left\{ \det 
		\Bigl(W(\la) \Bigr)\exp \Bigl(-\bigl[\log\det W\bigr](\la) \Bigr) 
		\right\}= 
		1 }$, where $W(\la)$ and $[\log\det W\bigr](\la) $ are \linebreak harmonic extensions of 
		functions $W\Bigm|  \T$ and $\log\det W \Bigm |\T$ respectively at point 
		$\la\in\D$. 
	
		\item  The measures
		$$
\left\| W(z)^{-1/2} \left(\frac{\partial}{\partial x} W(z)\right) W(z)^{-1/2}
\right\|^2(1-|z|^2) dxdy
$$ 
and 
$$
\left\| W(z)^{-1/2} \left(\frac{\partial}{\partial y} W(z)\right) 
W(z)^{-1/2}
\right\|^2(1-|z|^2) dxdy
$$ 
are vanishing Carleson measures
	\end{enumerate}
\end{thm}

Together with Theorem \ref{t-Pe1} the above theorem yields the 
complete description of completely regular stationary processes

\begin{thm}
	A stationary process with spectral density $\mathcal{W}$ is completely 
	regular if and only if $\mathcal{W}$ admits the representation
	$$
	\mathcal{W} = P^{*} W P,
	$$
	where $P$ is a polynomial matrix whose determinant has roots on 
	$\T$  and the matrix-function $W$ satisfies $W^{-1}\in L^{1}$ and 
	one of equivalent conditions 3--6 of Theorem \ref{t-main}
\end{thm}

Let us discuss the main result (Theorem \ref{t-main}) a little bit. 
First of all it is not difficult to show directly that in the scalar 
case the conditions 3--6 of Theorem \ref{t-main} are equivalent to 
$W=e^{\varphi}$, $\varphi\in \text{VMO}$. We are leaving this as an 
exercise for the reader. 

Usually in probability only real valued stationary processes are 
considered. In that case the spectral density of a process should 
satisfy $W(\overline z) = W(z)$, and only such functions can be 
realized as 
densities of stationary processes. 

If one allow complex-valued processes, any non-negative matrix 
function is the spectral density of some stationary process. 

Our theorem deals with arbitrary non-negative matrix-functions and can 
be applied to complex-valued processes (as well as to real-valued). 

\section{Scheme of the proof of the main result}

The diagram of the proof will be the following: $1\Longrightarrow 4  
\Longrightarrow 5 \Longrightarrow 6 \Longrightarrow 1 $. Then we will 
show that  $1\Longrightarrow 2$ and so automatically $2\Longrightarrow 
1$. 

And in this section we will show that $3\Longleftrightarrow 4$. 

\begin{lm}
\label{l2.1}
For a scalar weight $w$ the following conditions are equivalent: 
\begin{enumerate}
\item ${\displaystyle \limsup_{|I|\to 0}\Bigl(\frac{1}{|I|} \int_{I} w \Bigr) 
\Bigl(\frac{1}{|I|} \int_{I} w^{-1} \Bigr) = 1}$;
\item ${\displaystyle \limsup_{|\la|\to1} w(\la) w^{-1}(\la) = 1}$, 
where $w(\la)$ and $w^{-1}(\la)$ denote the harmonic extensions of 
$w$ and $w^{-1}$ respectively at the point $\la$; 
\item ${\displaystyle w=e^{\f}}$, where $\f\in\text{VMO}$. 
\end{enumerate}
\end{lm}

\begin{proof}
First of all let us rewrite  condition $1$. Let $\f:= \log w$.  
For a function $f$ let $f\ci I$ denote its average over the arc $I$,
$f\ci I := |I|^{-1}\int_I f$. Then clearly
$$
w\ci I\cdot (w^{-1})\ci I = \bigl[ w\ci I \exp(-\f\ci I) \bigr] \cdot
\bigl[(w^{-1})\ci I \exp(\f\ci I) \bigr]. 
$$
By Jensen inequality (geometric mean $\le$ arithmetic mean) the
expressions in brackets are at least $1$, so the condition ! splits
into the following 2 conditions
$$
\limsup_{|I|\to 0} \bigl[ w\ci I \exp(-\f\ci I) \bigr] = 1, \qquad
\text{and} \qquad \limsup_{|I|\to 0} \bigl[ w^{-1}\ci I \exp(\f\ci I)
\bigr] = 1. 
$$
Let $f_+$ denote the positive part of the function $f$, 
$f_+(x) := \max(f(x),0)$. Then the inequality 
$$
x\le e^x-1 \qquad \text{for} \ x\ge 0
$$
implies
$$
\frac{1}{|I|}\int_I (\f-\f\ci I)_+ \le \frac{1}{|I|}\int_I
\Bigl(  \exp(\f-\f\ci I) - 1 \Bigr) =
 w\ci I \exp(-\f\ci I) - 1 \to 0 \quad \text{as}\ |I|\to 0.  
$$
Since $\int_I |\f-\f\ci I | = 2 \int_I (f-f\ci I)_+$, one can conclude
that $\f\in \text{VMO}$. 

Similarly, using Poisson averages instead of averages over intervals 
one can get from condition 2 of the lemma that harmonic extension of $|\f-\f(\la)|$ at the point 
$\la$ tends to $0$ as $\la\to1$. But that is an equivalent definition 
of VMO, so the condition 2 also implies that $\f\in\text{VMO}$. 

On the other hand, if $\f\in\text{VMO}$, John--Nirenberg Theorem (see 
\cite{Gar}[Chapter VI], the measure of the set $\{t\in I: |\f(t)-\f\ci I|> 
a\}$ is estimated from above by $Ce^{-K a }$, where 
$K=K\ci I\to \infty$ as $|I|\to 0$. Therefore for $x>1$ the measure of 
the set $\{ t\in I: \exp(\f(t)-\f\ci I) >x \}$ is estimated from above 
by $C x^{-K}$. Integrating this distribution function one can get 
that  $\limsup_{|I|\to0}w\ci I \exp(-\f\ci I) \le 1$ (in fact, it is
$1$, because by Jensen inequality $w\ci I \exp(-\f\ci I)\ge 1$. 
Similarly, $ \limsup_{|I|\to0}(w^{-1})\ci I \exp(\f\ci I) = 1$. 
Multiplying the above two inequalities one gets condition 1. 

The proof that $3\Longrightarrow 2$ is similar. For a point 
$\la\in\D$ let $I_{\la}$ be an interval with center at $\la/|\la|$ of 
length $\sqrt{1-|\la|}$. Since the Poisson Kernel $P_{\la}(z) = 
(1-|\la|^{2})\cdot |1-\overline\la z|^{-2}$ satisfies 
$\sup_{z\in\T\setminus I_{\la}} P_{\la}(z) \to 0$ as $|\la|\to1$, the 
distribution inequality for $\f$ on $I_{\la}$ implies that 
$w(\la)\cdot\exp(-\f(\la)) \to 0$ as $|\la|\to 1$, and therefore the 
condition 2 of the lemma. 
\end{proof}

The following Lemma is probably well known and can be easily from the 
distribution function inequality for VMO (John--Nirenberg Theorem). 

\begin{lm}
\label{l2.2}For $\la\in\D$ let $I_{\la}$ be an interval centered at 
$\la/|\la|$ of length $1-|\la|$. If $\f\in\text{VMO}$, then 
$\f\ci{I_{\la}} - \f(\la) \to 0$ as $|\la|\to1$. 
\end{lm}

\begin{cor}
\label{c2.3}
Let $\f\in\text{VMO}$ and let $w=e^{\f}$. Then for $I_{\la}$ as in 
the above lemma we have 
$$
\lim_{|\la|\to1} \frac{w(\la)}{w\ci{I_{\la}}} =1. 
$$
\end{cor}
\begin{proof}
By the above lemma $\lim_{|\la|\to1} 
\exp(w(\la))/\exp(w\ci{I_{\la}}) =1$.  On the other hand it follows 
from the proof of Lemma \ref{l2.1} that 
$$
\lim_{|\la|\to1} w(\la)/\exp(\f(\la)) =1 \qquad \text{and} \qquad 
\lim_{|I|\to0} w\ci I /\exp(\f\ci I) =1. 
$$	
Taking the ration of the last 2 identities (with $I=I_{\la}$)we get the 
statement we need. 
\end{proof}

Now to show equivalence of condition 3 and 4 of the Theorem 
\ref{t-main} is  enough to show that these conditions imply that for 
a fixed vector $e\in\C^{d}$ scalar weight $w(z)=(W(z)e,e)$ satisfies 
conditions 1 and 2 of Lemma \ref{l2.1}. Then Corollary \ref{c2.3} 
implies that the averages $W\ci{I_{\la}}$ and $W(\la)$ are 
equivalent, the same holds for $W^{-1}$, and we are done. 

It remains now to show that the scalar weight 
$w(z)=(W(z)e,e)$ satisfies condition 1 (equivalently 2) of Lemma 
\ref{l2.1}. The easiest way to do that is to recall where the 
Muckenhoupt condition $(A_{2})$ came from, see \cite{TV1}. 

Recall that the quantity $\bigl\| [W\ci I]^{1/2}[(W^{-1})\ci 
I]^{1/2} \bigr\|$ is just the norm of the averaging operator $f\mapsto 
f\ci I\cdot \chi\ci I$ in the weighted space $L^{2}(W)$, see 
\cite{TV1}[Lemma 2.1]. Then $[w\ci 
I]^{1/2}[(w^{-1})\ci I]^{1/2}$ is the norm of the restriction of 
the above averaging operator onto the subspace of $L^{2}(W)$ 
consisting of functions of form $fe$ where $f$ is a scalar function. 
Therefore 
$$
1\le [w\ci 
I]^{1/2}[(w^{-1})\ci I]^{1/2} \le \bigl\| [W\ci I]^{1/2}[(W^{-1})\ci 
I]^{1/2} \bigr\|
$$
so the weight $w$ satisfies condition 1 of the lemma.

Similarly, the quantity 
$\bigl\| W(\la)^{1/2}W^{-1}(\la)^{1/2} \bigr\|$ is just the norm of 
another averaging operator $\bigl(f\mapsto \int_{\T} 
fk_{\la}\bigr)k_{\la}$, where $k_{\la}$ is the normalized reproducing 
kernel of $H^{2}$, $k_{\la}(z)=(1-|\la|^{2})^{1/2}(1-\overline \la 
z )^{-1}$, see \cite{TV1}[Lemma 2.1], so condition 4 of the theorem 
imply condition 2 of the lemma for the weight $w$.

\section{Eliminating probability}

The problem of description of completely regular processes can be now
stated without mentioning any probability theory at all.

First of all notice that without loss of generality we can assume that
the process is complex-valued. Namely, if we have a real stationary
process $\{x(n)\}_{n\in \Z}$ we can consider its comlexification,
namely the same process but in the complex Hilbert space
$L^2_{\C^d}(\Omega, dP)$. Considers the comlexificated
past $(\mathcal{X}_n)\ci{\C}$  and future $(\mathcal{X}^n)\ci{\C}$
\begin{align*}
	(\mathcal{X}_{n})_{\C} & = \spn\left\{ x_{j}(k)\,:\, 1\le j \le d,\, 
	k<n\right\}
	\\
	(\mathcal{X}^{n})_{\C} & = \spn\left\{ x_{j}(k)\,:\, 1\le j \le d,\, 
	k\ge n\right\}
\end{align*}
where $\spn$ now means the closed linear span in the complex Hilbert
space $L^2_{\C^d}(\Omega, dP)$.  It is easy to see that 
\begin{multline*}
\sup\left\{ | \bE (\xi {\eta}) |\,:\, \xi 
\in\mathcal{X}_{0},\, \eta\in \mathcal{X}^{n},\, \bE |\xi|^{2}\leq 
1,\, \bE |\eta|^{2}\leq 1 \right\} =
\\
= \sup\left\{ | \bE (\xi\bar{{\eta}}) |\,:\, \xi 
\in(\mathcal{X}_{0})_{\C},\, \eta\in (\mathcal{X}^{n})_{\C},\, \bE |\xi|^{2}\leq 
1,\, \bE |\eta|^{2}\leq 1 \right\}\,,
\end{multline*}
so a process and its comlexification are completely regular
simultaneously. So we indeed can assume from the beginning that our
process is complex valued.

Consider now the vector space $L^2(W)$ of $\C^d$-valued  functions on the
unit circle with the norm
$$
\| f \|_{L^2(W)}^2 = \int_\T (W(\xi) f(\xi), f(\xi))\ci{\C^{d}} dm(\xi)
$$
(of course we have to take the quotient space over the functions of
norm $0$). The mapping $x_j(k) \mapsto z^k e_j$, where $e_j$,
$j=1,...,d$ is the standard orthonormal basis in $\C^d$, is an isometric
isomorphism between $\spn\{x_j(k)\,:\, 1\le j \le d, k\in \Z\}$ and
$L^2(W)$. 

The past  $\mathcal X_n$ and future $\mathcal X^n$ are mapped to the spaces 
$X_n$ and $X^n$ of $L^2(W)$
\begin{align}
X_n &= \spn\{ z^k \C^d\,:\, k<n \}\, \label{eq-past}\\ 
X^n &= \spn\{ z^k \C^d\,:\, k\ge n \}\,. \label{eq-future}
\end{align}
 
So the problem  of describing completely regular stationary processes
can be reformulated as follows: {\em describe all matrix weights $W$ such
that the spaces $X_0$ and $X^n$ are asymptotically (as $n\to \infty$)
orthogonal to each other,} 
\begin{equation}
\label{def-rho}
\rho_n= \sup\left\{ |  (\xi, {\eta})\ci{L^2(W)} |\,:\, \xi 
\in{X}_{0},\, \eta\in {X}^{n},\,  \|\xi\|\ci{L^2(W)}\leq 
1,\, \|\eta\|\ci{L^2(W)}\leq 1 \right\} \longrightarrow 0\,,
\end{equation}
as $n\to \infty$.

\section{Necessity ($1\Longrightarrow 4$)}
\label{s-nec}
\setcounter{equation}{0}

In this section we are going to prove the implication $1\Longrightarrow 4$ 
(see Theorem \ref{t-nec} below) and the equivalence $1\Longleftrightarrow 
2$ (see Lemma \ref{l-W-W1}).

 For a
function $F$ defined on the unit circle $\T$ let  $F(\la)$ denote its harmonic
extension a the point $\la \in \D$. 

\begin{thm}
\label{t-nec}

Let $W$ be a matrix valued weight such that $W^{-1}\in L^1$. Suppose the
``past'' $X_0$ and ``future'' $X^n$ defined by \eqref{eq-past},
\eqref{eq-future} are asymptotically orthogonal, that means
$$
\rho_n= \sup\left\{ |  (\xi, {\eta})\ci{L^2(W)} |\,:\, \xi 
\in{X}_{0},\, \eta\in {X}^{n},\,  \|\xi\|\ci{L^2(W)}\leq 
1,\, \|\eta\|\ci{L^2(W)}\leq 1 \right\} \longrightarrow 0\
$$
as $n\to \infty$. Then 
$$
\limsup_{|\la|\to 1}\left\| 
		\Bigl(W(\la) \Bigr)^{1/2} \Bigl(W^{-1}(\la)
		\Bigr)^{1/2} 
\right\| = 1\,. 
$$
\end{thm}

\begin{proof}
First of all let us show that if $W^{-1}$ is completely regular and
$W^{-1}\in L^1$ then $W$ satisfies the Muckenhoupt ($A_2$)  condition
\begin{equation}
\sup_{\la\in\D} \left\| 
		\Bigl(W(\la) \Bigr)^{1/2} \Bigl(W^{-1}(\la)
		\Bigr)^{1/2} 
\right\| < \infty \,. 
\tag{$A_p$}
\end{equation}
 Recall that $\left\| 
		\Bigl(W(\la) \Bigr)^{1/2} \Bigl(W^{-1}(\la)
		\Bigr)^{1/2} 
\right\|$ is exactly the norm of the operator $f\mapsto
(f,k_\la)k_\la$ in the weighted space $L^2(W)$; here
$k_\la$ denotes the normalized  reproducing kernel for $H^2$, 
$$
k_\la(z) := \frac{(1-|\la|^2)^{1/2}}{1-\overline\la z}\,, \qquad
\la\in \D,
$$
$\|k_\la\|_2=1$. Note that $k_0\equiv1$. So if $W^{-1}\in L^1$ the
operator
$f\mapsto (f,1)1$ is bounded in $L^2(W)$, and therefore by translation invariance
the operators $f\mapsto (f,z^n)z^n=\hat f(n)z^n$ are bounded as well (they all have
the same norm).

We know that the spaces $X_0$ and $X^n$ are asymptotically orthogonal,
so we can say that for large enough $N$ the operator $P_+$ restricted
onto $\spn\{X_0,X^N\}= \spn\{ z^n\C^d\,:\,n\notin [0,N]\}$ is bounded,
say by $2$, 
$$
\|P_+ f\|\ci{L^2(W)} \le 2\|f\|\ci{L^2(W)}\,,\qquad \forall f\in \spn\{X_0,X^N\}= \spn\{ z^n\C^d\,:\,n\notin [0,N]\}\,.
$$
Since $f-\sum_{n=0}^N \hat f(n)z^n \in \spn\{X_0,X^N\}= \spn\{
z^n\C^d\,:\,n\notin [0,N]\}$, one can conclude that the operator $P_+$
is bounded in $L^2(W)$ and so the weight satisfy the Muckenhoupt
condition $(A_2)$.

We will need the following simple lemma about Muckenhoupt weights. 

\begin{lm}
\label{l-doubl}
If $w$ is a scalar Muckenhoupt weight, then its harmonic extension
$w(\la)$ can't decay to fast near the boundary of the disk.  Namely,
if the Muckenhoupt norm of $w$ is at most $M$ there is a function
$\alpha=\alpha\ci{M}$, 
$\alpha:[0,1)\to (0,\infty)$, $\alpha(t)\searrow 0$ as $t\to 1+$ such that
$$
\frac{(1-|\la|^2)w(0)}{w(\la)} \le \alpha(|\la|).
$$
\end{lm}

\begin{proof}[Proof of the lemma]
For an arc $I\subset \T$ and $k>0$ let $kI$ denote the arc of length
$k|I|$ with the
same center as $I$. 

We are going to show that for a Muckenhoupt weight $w$ with the
Muckenhoupt norm at most $M$ 
\begin{equation}
\label{eq-doubl}
w\ci{2^n I} \le M^2 (2-\varepsilon)^n w\ci{I}\,, \qquad \e=\e(M)>0\,.
\end{equation}
Applying this formula in the case $2^n I =\T$ and using the trivial estimate
$$
w(\la) \ge C w\ci{I_\la}
$$
where $I_\la$ is the arc with center at the point $\la/|\la|$,
$|I_\la|=1-|\la|^2$ and $C$ is an absolute constant, we can get 
from there (recall that $|I_\la|=1-|\la|^2=2^{-n}$)
$$
w(\la) \ge c(2-\e)^{-n}\cdot w(0) = 
c (2-\e)^{\log_{2}(1-|\la|^2)} \cdot w(0) = 
 c\cdot(e-\delta)^{\log(1-|\la|^2)}\cdot w(0),  
$$
where $\delta=\delta(\e)>0$; here $e$ is the base of the natural 
logarithm, not a vector in $\C^{d}$. 
This estimate implies the conclusion of the lemma with
$\alpha(t)=c^{-1}(1-t^2) \cdot (e-\delta)^{-\log(1-t^2)}$.

To prove \eqref{eq-doubl} we notice the since the weight $w^{-1}$ is
the Muckenhoupt $(A_2)$ weight with the same Muckenhoupt norm as $w$,
it is doubling and therefore 
$$
(w^{-1})\ci{2I} \ge (2-\e)^{-1} (w^{-1})\ci{I}\,,
$$
where $\e$ depends only on the Muckenhoupt norm of $w$. Iterating this
inequality $n$ times we get 
$$
(w^{-1})\ci{2^n I} \ge (2-\e)^{-n} (w^{-1})\ci{I}\,.
$$
The last estimate and the Muckenhoupt condition imply
$$
w\ci{2^n I} \le M/(w^{-1})\ci{2^n I} \le M\cdot(2-\e)^n/(w^{-1})\ci{I}
\le M^2 w\ci{I}
$$
and that is exactly what we need.
\end{proof}

\begin{cor}
\label{c-doubl}
If a matrix weight $W$ satisfies the Muckenhoupt condition $(A_2)$
with the Muckenhoupt norm at most $M$
then for any $e\in\C^d$
$$
(1-|\la|^{2})\cdot \frac{(W(0)e,e)\ci{\C^d}}{(W(\la)e,e)\ci{\C^d}} \le \alpha(|\la|) \
\to 0 \qquad \text{as}\ |\la|\to 1, 
$$
where $\alpha=\alpha_M$ is the function from Lemma \ref{l-doubl}.
\end{cor}
\begin{proof}[Proof of the corollary] The proof follows immediately
from the fact that the scalar weight $w$, $w(\xi) =
\bigl(W(\xi)e,e\bigr)_{\C^d}$ is the Muckenhoupt $(A_2)$ weight with
the Muckenhoupt norm at most $M$ (see \cite{TV2}, proof of Corollary
2.4).  
\end{proof}

We now return to the proof of the theorem.

The condition $W^{-1}\in L^1$ implies that $\int_\T\log\det W(\xi)
dm(\xi)>-\infty$, hence (see \cite{Ro}) 
there exists a factorization
of $W$ of the form $W=F^*F$, where $F$ is an outer matrix function in
$H^2$. 

Take $e\in\C^d$ and let us compute the distance
$$
\dist\ci{L^2(W)}\{z^{-1}e, \spn\{z^n\C^d\,:\,n\ge 0\} =
\dist\ci{L^2(W)}\{e, \spn\{z^n\C^d\,:\,n> 0\}\,.
$$
By the vectorial version of the Szeg\"{o} theorem (see \cite{Ro}) 
this
distance is exactly $\|F(0)e\|$. Using the M\"{o}bius transformation of
the disk one can get from there
$$
\dist\ci{L^2(W)}\{ \frac{(1-|\la|^2)^{1/2}}{z-\la}e,
\spn\{z^n\C^d\,:\,n\ge 0\}\} = \|F(\la)e\|\ci{\C^d}\,.
$$

Writing the Fourier series expansion of
$\frac{(1-|\la|^2)^{1/2}}{z-\la}$ 
$$
\frac{(1-|\la|^2)^{1/2}}{z-\la}= (1-|\la|^2)^{1/2}\sum_{n=0}^\infty
\la^n z^{-(n+1)}
$$
one can see that for any fixed $N>0$ the function
$\frac{(1-|\la|^2)^{1/2}}{z-\la} e$ is almost in the ``past''
$X\ci{-N}$ as $|\la|\to 1$. Namely, 
$$
f_\la = \frac{(1-|\la|^2)^{1/2}}{z-\la} e = (1-|\la|^2)^{1/2}\sum_{n=0}^{N-1}
\la^n z^{-(n+1)}e + (1-|\la|^2)^{1/2}\sum_{n=N}^{\infty}
\la^n z^{-(n+1)}e=f_\la^1+f_\la^2\,,
$$
where $f_\la^2\in X\ci{-N}$, and $f_\la^1$ is small, 
$$
\frac{\|f_\la^1\|\ci{L^2(W)}}{\| f_\la\|\ci{L^2(W)}} 
 \le \frac{ (1-|\la|^2)^{1/2} N\cdot 
\|e\|\ci{L^2(W)}}{ \bigl(W(\la)e,e\bigr)\ci{\C^d}^{1/2}} = 
\frac{(1-|\la|^2)^{1/2} N\cdot
 \bigl(W(0)e,e\bigr)\ci{\C^d}^{1/2}}{\bigl(W(\la)e,e\bigr)\ci{\C^d}^{1/2}} 
\le N \alpha(|\la|)^{1/2}\to 0,
$$ 
as $|\la|\to 1$, where $\alpha(.)$ is as in Lemma \ref{l-doubl} and Corollary
\ref{c-doubl}. 

Since $X_{0}$ and $X^{N}$ are asymptotically orthogonal, the shift invariance 
implies that the  subspaces 
$X_{-N}$ and $X^{0}$ are asymptotically orthogonal as well. Taking 
$|\la|\to1$ and then $N\to\infty$ we can 
conclude that 
\begin{multline*}
\|F(\la)e\|\ci{\C^{d}} /\|W(\la)^{1/2}e\|\ci{\C^d} = \\=
\dist\ci{L^{2}(W)} \{ \frac{(1-|\la|^2)^{1/2}}{z-\la}e,
\spn\{z^n\C^d\,:\,n\ge 0\}\} / \left\| f_\la  
\right\|\ci{L^{2}(W)} \ge 1- \beta(|\la|)^{1/2} \to 1\,,
\end{multline*}
where $\beta(.)$ depends only on the Muckenhoupt norm  of $W$ and
$\beta(|\la|)\to 0$ as $ |\la|\to1$. 

The last inequality   implies 
\begin{equation}
\label{WF}
\|W(\la)^{1/2}F(\la)^{-1}\|\le (1-\beta(|\la|))^{-1}\,.
\end{equation}

Note that since $\|F(\la)e\|\ci{\C^{d}} /\|W(\la)^{1/2}e\|\ci{\C^d} 
\le 1$
for all $e\in\C^d$, we have 
$$
\|W(\la)^{1/2}F(\la)^{-1}\|\ge 1\,. 
$$ 

We will show a little later that under assumptions of the theorem the 
subspaces $X_{0}$ and $X^{N}$ in the weighted space $L^2(W^{-1})$ are
asymptotically orthogonal as well. The factorization $W=F^*F$ yields
the factorization $W^{-1}=F^{-1}(F^{-1})^*$ of $W^{-1}$. Similarly to
the previous case
$$
\dist\ci{L^2(W^{-1})} \{ \frac{(1-|\la|^2)^{1/2}}{1-\overline\la z}e,
\spn\{z^n\C^d\,:\,n\ge 0\}\} = \|F^{-1}(\la)^*e\|\ci{\C^d}
= \|F(\la)^{-1*}e\|\ci{\C^d}. 
$$
Acting as before we get
\begin{equation}
\label{W-1-F}
\|W^{-1}(\la)^{1/2}F(\la)^{*}\|\le (1-\beta_1(|\la|))^{-1}\,
\end{equation}
where $\beta_1(|\la|)\to 0$ as $ |\la|\to1$. 

Combining \eqref{WF} and \eqref{W-1-F} we get 
$$
\| W(\la)^{1/2}W^{-1}(\la)^{1/2} \| \le 
(1-\beta(|\la|))^{-1}(1-\beta_{1}(|\la|))^{-1} \to 1 \qquad 
\text{as} \ |\la|\to 1\,.
$$
So, we completed the proof modulo the following lemma. 
\end{proof}

This lemma also gives us the equivalence $1\Longleftrightarrow 2$. 
\begin{lm}
\label{l-W-W1}
Under assumptions of Theorem \ref{t-nec} the weight $W^{-1}$ is a 
spectral density of a completely regular process, i.e that the 
spaces $X_{0}$ and $X^{N}$ are asymptotically orthogonal (as 
$N\to\infty$) in the weighted space $L^{2}(W^{-1})$.
\end{lm} 

\begin{proof}

It is enough to show that 
$$
\| P_{+}\bigm| \spn\{X_{0},X^{N}\}\|\ci{L^{2}(W^{-1})\to L^{2}(W^{-1})} 
\to 1\qquad \text{as}\ N\to\infty.
$$
The later is true because
\begin{multline*}
\| P_{+}\bigm| \spn\{X_{0},X^{N}\}\|\ci{L^{2}(W^{-1})\to 
L^{2}(W^{-1})} =
\|W^{-1/2} \bigl( P_{+}\bigm| \spn\{X_{0},X^{N}\} \bigr) 
W^{1/2}\|\ci{L^{2}\to L^{2}} 
	\\ =
\|W^{1/2} \bigl( P_{+}\bigm| \spn\{X_{0},X^{N}\} \bigr) 
W^{-1/2}\|\ci{L^{2}\to L^{2}} 	
=
\| P_{+}\bigm| \spn\{X_{0},X^{N}\}\|\ci{L^{2}(W)\to 
L^{2}(W)}
\end{multline*}
and 
$$
\| P_{+}\bigm| \spn\{X_{0},X^{N}\}\|\ci{L^{2}(W^{-1})\to 
L^{2}(W)} \to 1 \qquad \text{as}\ N\to\infty
$$
(since $X_{0}$ and $X^{N}$ are asymptotically orthogonal in 
$L^{2}(W)$). 
\end{proof}

\section{Vanishing Carleson measures}
\setcounter{equation}{0}

Recall that $W(\la)$ and $W^{-1}(\la)$ denote harmonic extensions at 
the point $\la\in\D$ of the weights $W$ and $W^{-1}$ respectively.
 
\begin{lm}
\label{l-A-infty}
Let a matrix weight $W$ satisfy 
$$
\lim_{|\la|\to 1} \| W(\la)^{1/2}\bigl(W^{-1}\bigr)(\la)^{1/2} \| =1.
$$
Then 
$$
\limsup_{|\la|\to 1}\left\{ \det \Bigl(W(\la) \Bigr)\exp 
\Bigl(-\bigl[\log\det W\bigr](\la) \Bigr) \right\}= 1\,.
$$
\end{lm}
\begin{proof}
First of all let us notice that the assumption of the lemma implies that 
$W, W^{-1}\in
L^1(\T)$, therefore $\log(\det W) \in L^1(\T)$. Therefore there 
exists a
factorization $W=F^*F$ a.e. on $\T$, where $F$ is an outer function in
$H^2(M_{d\times d})$ 

Since $F$ is an outer function in $H^2$, $\det F$ is an outer 
function in
$H^{2/d}$. Therefore
\begin{equation}
\label{2.06}
|\det F(z)|= \exp \left\{ \left(\log|\det F| \right)(z) \right\} =
\exp\left\{\frac12\left(\log\det W \right)(z) \right\}
\end{equation}

It is well known fact that $F^*(z)F(z)\le W(z)$ for any $z\in D$, 
where $\le $
means the inequality for quadratic forms.  There are many proofs of 
this fact,
for example it admits a very simple operator-theoretic interpretation 
which is in fact hidden in the proof of Theorem \ref{t-nec}.
Explanation that we present here is more function-theoretic: Direct 
computation
shows that
$$
\Delta \left( F(z)^*F(z) \right) = 4\left(\bar\partial F(z)^*\right) 
\left (\partial F(z)\right) = 4 \left(\partial F(z)\right)^* 
\left (\partial F(z)\right) \ge 0\,,
$$
so for any $e\in \C^d$ the function $\|F(z)e\|^2$ is subharmonic and 
coincide
with $(W(\xi)e,e)$ on $\T$. 

We can do the same factorization for $W^{-1}$.  Namely, let $G$ be an 
outer matrix-valued function in $H^2(M_{d\times d})$ such that 
$W^{-1}=G^*G$ on $\T$.  We should point out to the reader that in 
general $G$ does not necessarily coincide with $F^{-1}$.  However, 
applying (\ref{2.06}) to $G$ one can conclude that
\begin{equation}
\label{2.07}
|\det G(z)|=\exp\left\{\frac12\left(\log\det W^{-1} \right)(z) 
\right\}
=|\det F(z)|^{-1}
\end{equation}
Now we are in position to prove the lemma. By the assumption
\begin{equation}
\lim_{|z|\to 1} \left\|  W(z)^{1/2}( W^{-1})(z)^{1/2} \right\| = 1,
\end{equation}
and therefore, 
$$
\lim_{|z|\to 1} \left|\det(W(z))\det\left((W^{-1})(z)\right)\right| 
= 1
$$  
Using (\ref{2.07}) one can rewrite the last identity as 
$$
\lim_{|z|\to 1} 
\left\{
  \left[ \det W(z) / |\det F(z)|^2 \right]
  \left[ \det W^{-1}(z) / |\det G(z)|^2 \right]
\right\} =1
$$
Since $F(z)^*F(z) \le W(z)$ and $G(z)^*G(z)\le W^{-1}(z)$, 
expressions in
brackets are at least 1, so, taking into account (\ref{2.06}) 
$$
\lim_{|z|\to 1} 
  \left[  \det W(z) / \exp \left\{ (\log \det W)(z) \right\} %
  \right]\ = 0
$$
or equivalently
\begin{equation}
\label{2.1}
\lim_{|z|\to 1} \log \left\{  \det( W(z) )\right\} - 
     \left( \log \det W \right) (z) =0\,.
\end{equation}

\end{proof}


\begin{thm}
\label{VCM}
A matrix weight $W$ satisfy 
$$
\limsup_{|\la|\to 1}\left\{ \det \Bigl(W(\la) \Bigr)\exp 
\Bigl(-\bigl[\log\det W\bigr](\la) \Bigr) \right\}= 1
$$
if and only if  the measures 
$$
\left\| W(z)^{-1/2} \left( \frac{\partial}{\partial x}W(z)\right) W(z)^{-1/2}
\right\|^2(1-|z|^2) dxdy
$$
and 
$$
\left\| W(z)^{-1/2} \left( \frac{\partial}{\partial y}W(z)\right) W(z)^{-1/2}
\right\|^2(1-|z|^2) dxdy
$$
are vanishing Carleson measures. 
\end{thm}

The implication $3 \Longrightarrow 4$ of Theorem \ref{t-main} follows immediately 
from Theorem \ref{VCM} and Lemma \ref{l-A-infty}. 

To prove the theorem we need the following well known description of 
vanishing Carleson measures
\begin{lm}
\label{l-VCM}
A measure $\mu$ in the unit disk $\D$ is a vanishing Carleson measure 
if and only if 
$$
\lim_{|\la|\to 1} \int_{\D} \frac{1-|\la|^{2}}{|1-\overline\la 
z|^{2}} d\mu(z) = 0.
$$
\end{lm}

We also need the following lemma that was proved in \cite{TV2}, see 
Lemma 3.1 there. 
\begin{lm}
\label{l2.1-TV2}
Let $W$ be a harmonic function of $n$ variables with values in the 
space 
of strictly positive $d\times d$ matrices ($W(x)=W(x)^*>0$ $\forall 
x$). Then 
$$
\Delta \left( \log(\det W) \right) = -\sum_{j=1}^n\trace\left(
(W^{-1/2}\frac{\partial W}{\partial x_j} W^{-1/2})^2\right)
$$
\end{lm}

\begin{proof}[Proof of Theorem \ref{VCM}]
The proof below follows the lines of the proof of Theorem 3.2 of 
\cite{TV2}.

By Green's formula and Lemma \ref{l2.1-TV2} 
\begin{multline*}
\log \left\{ \det( W(s) ) \right\} - 
     \left( \log \det W \right) (s)  =
-\frac{1}{2\pi}\iint_\D \log \left| \frac{1- \overline{s}
z}{z-s}
\right|
          \Delta \log\left\{  \det(W(z)) \right\}\, dx dy =
\\
=
\frac{1}{4\pi}\iint_\D 
\left\{ \trace
\left(W(z)^{-1/2}\frac{\partial W(z)}{\partial x} 
W(z)^{-1/2}\right)^2 +
\right.
\\ +
\left. \trace
\left(W(z)^{-1/2}\frac{\partial W(z)}{\partial y} W(z)^{-1/2}\right)^2
\right\}
\log \left| \frac{1- \overline{s} z}{z-s} \right|^2 \, dx dy
\end{multline*}

Using an elementary inequality $\log(1/a) \ge 1-a$ for $0<a\le 1$ and 
the fact
that $\|A\|\le\trace A$ for a non-negative matrix $A$,  
 the last integral is at least 
\begin{multline*}
\frac{1}{4\pi}\iint_\D 
\left\|
   W(z)^{-1/2}\frac{\partial W(z)}{\partial x} W(z)^{-1/2}
\right\|^2
\log \left| \frac{1- \overline{s} z}{z-s} \right|^2 \, dx dy
\ge \\
\frac1{4\pi} \iint_\D  
\left\|
   W(z)^{-1/2}\frac{\partial W(z)}{\partial x} W(z)^{-1/2}
\right\|^2
\left(1- \left| \frac{1- \overline{s} z}{z-s} \right|^2 \right) \, dx 
dy
= \\
=
\iint_\D 
\left\|
   W(z)^{-1/2}\frac{\partial W(z)}{\partial x} W(z)^{-1/2}
\right\|^2 
\cdot 
\frac{(1-|s|^2)(1-|z|^2)}{|1-\overline{s} z|^2} \, dx dy
\end{multline*}
Together with (\ref{2.1}) this  imply 
$$
\lim_{|s|\to 1}\iint_\D   
\frac{(1-|s|^2)}{|1-\overline{s} z|^2 } \cdot 
\left\|
   W(z)^{-1/2}\frac{\partial W(z)}{\partial x} W(z)^{-1/2}
\right\|^2
(1-|z|^2)\, dx dy = 0
$$
that yields that the measure $\left\|
   W(z)^{-1/2}\left(\frac{\partial}{\partial x} W(z) \right) 
W(z)^{-1/2}
\right\|^2(1-|z|^2)\, dx dy$ is a vanishing Carleson measure.

The measure $\left\|
   W(z)^{-1/2}\left(\frac{\partial}{\partial y} W(z) \right) 
W(z)^{-1/2}
\right\|^2(1-|z|^2)\, dx dy$  is treated similarly.   

To prove the opposite implication, let us estimate the integral 
$$
\iint_\D 
 \trace
\left(W(z)^{-1/2}\frac{\partial W(z)}{\partial x} 
W(z)^{-1/2}\right)^2 
\log \left| \frac{1- \overline{s} z}{z-s} \right|^2 \, dx dy
$$
(the  integral with $\partial W/\partial y$ can be estimated absolutely the same way). Denote by $b_{s}$ a Blaschke factor with zero at the point $s$, 
$b_{s}(z)= (z-s)(1-\overline s z)^{-1}$. 

First of all, we can estimate the trace by $d\cdot \|\cdot\|$, where $d$ 
is dimension of the space. 
So we can estimate the integral by
$$
C\iint\limits_{\D} \left\| W(z)^{-1/2}\frac{\partial W(z)}{\partial x} 
W(z)^{-1/2}\right\|^2 \log |b_{s}(z)|^{-2} dxdy = 
\iint\limits_{|b_{s}(z)|<\e}\ldots \ + \iint\limits_{|b_{s}(z)|\ge \e}\ldots 
$$
To estimate the second integral we notice that 
$$
\log |b_{s}(z)|^{-2} dxdy \le C(\e) \frac{(1-|s|^2)(1-|z|^2)}{|1-\overline{s} z|^2}
$$
for $|b_{s}(z)|\ge\e$, and since the measure is a vanishing Carleson 
measure we can make the integral as small as we want when $|s|\to 1$. 

To estimate the first integral let make a trivial observation: if 
$w\in L^{1}(\T)$, $w\ge0$ and $w(z)$ denotes its harmonic extension 
at the point $z$, then for all $z$ such that $|z|\le 1/2$ (and 
therefore for all $z$ such that $|z|<\e\le 1/2$)
$$
\frac{\partial}{\partial x} w(z) \le C w(0)\,,
$$
where $C$ is an absolute constant. Combining this observation with the 
Harnack inequality $w(0)\le C'w(z)$, $|z|\le 1/2$, and applying it  to 
functions $w(.) =\bigl(W(\cdot)e,e)_{\C^{d}}$ we get the inequality 
for quadratic forms
$$
\frac{\partial}{\partial x}  W(z) \le C\e W(0) \le C_{1} W(z)\,.
$$
It in turn implies 
$$
\left\| W(z)^{-1/2} \Bigl(\frac{\partial}{\partial x}  W(z)\Bigr) W(z)^{-1/2}
\right\| \le C_{1}\,, \qquad \forall z:\ |z|< \e \le1/2\,.
$$
Using the M\"{o}bius transformation $z\mapsto b_{s}(z)$ we get 
$$
\left\| W(z)^{-1/2} \Bigl(\frac{\partial}{\partial x}  W(z)\Bigr) 
W(z)^{-1/2}
\right\| \le C\e\,, \qquad \forall z:\  |b_{s}(z)|< \e \le1/2\,.
$$
Since 
$$
\iint_{|b_{s}(z)|\le \e} \log |b_{s}(z)|^{-2} dxdy \le 
C\e^{2} \log\frac1\e\,,
$$
we can estimate the first integral by $C\e^{2} \log(1/\e)$; we can make 
this number as small as we want by picking sufficiently small $\e$.

\end{proof}

\section{Embedding theorem and equivalent norms}
\label{s-ImbThm}

By analogy with the scalar case (see \cite{TVZ}) we will say that a 
matrix weight $W$ satisfies the {\em invariant $A_{\infty}$ condition}
if 
\begin{equation}
\sup_{s\in\D} \left\{ \det \Bigl(W(s) \Bigr)\exp 
\Bigl(-\bigl[\log\det W\bigr](s) \Bigr) \right\}<\infty\,.
\tag{$\mathit{invA}_{\infty}$}
\end{equation}
The {\em supremum} is called the invariant $A_{\infty}$ norm of $W$. 

Theorem \ref{VCM} implies that if the measures
$$
\left\| W(z)^{-1/2} \left( \frac{\partial}{\partial x}W(z)\right) W(z)^{-1/2}
\right\|^2(1-|z|^2) dxdy
$$
and 
$$
\left\| W(z)^{-1/2} \left( \frac{\partial}{\partial y}W(z)\right) W(z)^{-1/2}
\right\|^2(1-|z|^2) dxdy
$$
are vanishing Carleson measures then the weight $W$ satisfies the 
invariant $A_{\infty}$ condition. 

Literally repeating the proof of Theorem \ref{VCM} one can obtain that 
the weight $W$ satisfies the invariant $A_{\infty}$ condition if and 
only if the above measures are Carleson.

We will need the following ``embedding theorem''. More general result 
was proved  in
\cite{TV2}, Lemma 4.1. 

\begin{lm}
\label{l3.1}
Let $W$ be a matrix  weight satisfying the invariant $A_{\infty}$ 
condition, and let  $\mu$ be a Carleson  measure with the Carleson
norm $\|\mu\|\ci{C}$.  Then for any analytic (or antianalytic) 
vector-function
$f$,  the  following inequality holds, 
$$
\iint_\D (W(z)f(z),f(z)) \,d\mu(z) \le C \|\mu\|\ci{C}\int_\T (W(\xi)f(\xi),f(\xi))
dm(\xi)\,,
$$
where the constant $C$ depends the dimension $d$ and  the invariant 
$A_{\infty}$ norm of $W$. 
\end{lm}

\begin{proof}
The invariant $A_{\infty}$ condition implies that $\log\det W\in L^{1}$, 
so there exists (see \cite{Ro}) an outer function $F\in H^{2}(M_{d\times d})$ such 
that $W=F^{*}F$. It is well known (see again \cite{Ro}) that 
$$
|\det F(z)| = \exp\left\{ \frac12 \bigl[ \log \det W \bigr](z) 
\right\}\,.
$$
It is well known and it was already shown it in the proof of Lemma 
\ref{l-A-infty} that $F(z)^{*}F(z) \le W(z)$. Hence 
\begin{equation}
\label{expand}
\| W(z)^{1/2}F(z)^{-1} e \| \ge \|e\|, \qquad e\in \C^{d}\,. 
\end{equation}
Since 
$$
\left|\det\left\{ W(z)^{1/2}F(z)^{-1} \right\}\right| = \left\{ \det \Bigl(W(\la) \Bigr)\exp 
\Bigl(-\bigl[\log\det W\bigr](\la) \Bigr) \right\}^{1/2} \le C
$$
we can estimate 
$$
\| W(z)^{1/2}F(z)^{-1} e \| \le C\,.
$$
Together with \eqref{expand} it implies that $(W(z)e,e)$ and 
$\|F(z)e\|^{2}$ are equivalent in a sense of two-sided estimate. 
Therefore
	\begin{multline*}
		\iint_\D (W(z)f(z),f(z)) \,d\mu(z) 
		\le C \iint_\D (F(z)f(z),F(z)f(z)) \,d\mu(z) \le \\
		\le C \|\mu\|\ci{C}\int_\T (F(\xi)f(\xi),F(\xi)f(\xi))dm(\xi) =
		C \|\mu\|\ci{C}\int_\T (W(\xi)f(\xi),f(\xi)) dm(\xi)\,.
	\end{multline*}	
\end{proof}

We also need the following simple lemma.

\begin{lm}[equivalence of weighted norms]
\label{t-norm-equiv}
Let $W$ be a matrix  weight satisfying the invariant $A_{\infty}$ 
condition. There exist a 
constant $C$ such that for any analytic or antianalytic
vector-function $f$ in
$L^2(W)$ satisfying $f(0)=0$
$$
\frac1C \int_\T (Wf,f)dm  \le 
\iint_\D 
 (W(z) f'(z), f'(z))  \log\frac{1}{|z|}\, dxdy 
\le 
C \int_\T (Wf,f)dm 
$$
\end{lm}

\begin{proof}Let us recall the the operators $\partial$ and $\overline\partial$ 
are defined as 
$$
\partial f = \frac12\left(
\frac{\partial f}{\partial x} - i \frac{\partial f}{\partial 
y}\right), \qquad 
\overline\partial f = \frac12\left(
\frac{\partial f}{\partial x} + i \frac{\partial f}{\partial 
y}\right)\,. 
$$
Recall  that for analytic functions $\partial f = f'$  and 
$\overline\partial f = 0$. 

 Let $f$ be an analytic function, $f(0)=0$. Using the 
Green's formula and taking into  account that $f(0)=0$ and $\Delta= 
4\partial \overline\partial = 4 \overline\partial \partial $ we get
\begin{multline*}
	\int_{\T}\bigl(Wf,f\bigr) dm =\frac{1}{2\pi} \iint_{\D} \Delta 
	\bigl(W(z)f(z),f(z)\bigr) \log\frac{1}{|z|} dxdy = \\
	\frac{2}{\pi}\iint_{\D}\bigl(\overline\partial W(z)f'(z),f(z)\bigr) \log\frac{1}{|z|} dxdy
+ \frac{2}{\pi}\iint_{\D}\bigl(\partial W(z)f(z),f'(z)\bigr) \log\frac{1}{|z|} dxdy
+\\
+ \frac{2}{\pi}\iint_{\D}\bigl( W(z)f'(z),f'(z)\bigr) \log\frac{1}{|z|} dxdy
= \frac2{\pi}(\mathcal I_{1} + \mathcal I_{2} + \mathcal I_{3})
\end{multline*}
The last integral $\mathcal I_{3}$ is exactly the integral we want to 
estimate. Let us denote $A^{2}:= \int_{\T}\bigl(Wf,f\bigr) dm $, 
$B^{2}:= \mathcal I_{3}$. We want to show that $A\asymp B$ in a sense 
of two sided estimate. Let us estimate $\mathcal I_{1}$:
\begin{multline*}
	| \mathcal I_{1} | = \left| \iint_{\D}
	\bigl(W(z)^{-1/2}\overline\partial 
	W(z)W(z)^{-1/2}W(z)^{1/2}f'(z),W(z)^{1/2}f(z)\bigr) \log\frac{1}{|z|} 
	dxdy \right|
	 \\
	\le \left| \iint_{\D}
	\bigl\| W(z)^{-1/2}\overline\partial 
	W(z)W(z)^{-1/2}\bigr\| \cdot \bigl\|W(z)^{1/2}f'(z)\bigr\| \cdot 
	\bigl\| W(z)^{1/2}f(z)\bigr\|\cdot  \log\frac{1}{|z|} dxdy \right| \\
	\le 
	\left( \iint_{\D}\bigl\| W(z)^{-1/2}\overline\partial 
	W(z)W(z)^{-1/2}\bigr\|^{2} \bigl( W(z) f(z), 
f(z) \bigr)_{\C^{d}} \log\frac{1}{|z|} dxdy \right)^{1/2}\times  \\
\times \left( \iint_{\D} \bigl(W(z) f'(z),f'(z)\bigr)_{\C^{d}} 
\log\frac{1}{|z|} dxdy \right)^{1/2}\,.
\end{multline*}
The measure $\bigl\| W(z)^{-1/2}\overline\partial W(z)W(z)^{-1/2}\bigr\|^{2}
\log\frac{1}{|z|} dxdy$ is Carleson, so by Lemma \ref{l3.1} the first 
term in the product is estimated by $KA$ ($K$ is a constant). The second term is just $B$ 
so 
$|\mathcal I_{1}| \le KAB$. Similarly $|\mathcal I_{2}| \le KAB$. So
$$
A^{2} = B^{2} + \mathcal I_{1} + \mathcal I_{2}\,,
$$
where 
$$
|\mathcal I_{1}|,\, |\mathcal I_{2}| \le KAB\,.
$$
This immediately implies 
$$
\frac1C A \le B \le CA
$$
for an appropriate choice of $C$. 
\end{proof}

\section{Proof of the implication $6\Longrightarrow 1$}

To prove the implication  $6\Longrightarrow 1$ we need to estimate 
$\int_{\T}(W z^{n}f,g)dm$, $f\in X^{0}$, $g\in X_{0}$, 
$\|f\|\ci{L^{2}(W)} = \|g\|\ci{L^{2}(W)} = 1$.

Using the 
Green's formula and taking into  account that $g(0)=0$ and $\Delta= 
4\partial \overline\partial = 4 \overline\partial \partial $ we get
\begin{multline*}
\int_{\T}(W z^{n}f,g)dm =\frac{1}{2\pi} \iint_{\D} \Delta \bigl(W(z) 
z^{n}f(z),g(z)\bigr)_{\C^{d}} \log\frac{1}{|z|} dxdy = \\
= \frac{2}{\pi}\iint_{\D} \bigl(\overline\partial W(z) 
\partial(z^{n}f(z)),g(z)\bigr)_{\C^{d}}\log\frac{1}{|z|} dxdy + \\ +
\frac{2}{\pi}\iint_{\D} \bigl(\overline\partial W(z) 
(z^{n}f(z)),\overline\partial g(z)\bigr)_{\C^{d}} \log\frac{1}{|z|} dxdy 
=\frac{2}{\pi}(\mathcal I_{1} + \mathcal I_{2})
\end{multline*}
The second integral is easy to estimate:
\begin{multline*}
|\mathcal I_{2} | \\
= \left|\iint_{\D} \bigl(W(z)^{-1/2}\overline\partial W(z) 
W(z)^{-1/2} W(z)^{1/2}
(z^{n}f(z)),W(z)^{1/2}\overline\partial g(z)\bigr)_{\C^{d}} \log\frac{1}{|z|} 
dxdy \right| \\
\le \iint_{\D} \| W(z)^{-1/2}\overline\partial W(z) W(z)^{-1/2}\| 
\cdot 
\| W(z)^{1/2} (z^{n}f(z)) \| \cdot  \|  W(z)^{1/2} \overline\partial g(z) \| 
\log\frac{1}{|z|} dxdy \\
\le \left(
\iint_{\D} |z|^{2n}\cdot \| W(z)^{-1/2}\overline\partial W(z) 
W(z)^{-1/2}\|^{2}  \cdot 
\bigl( W(z) f(z), 
f(z) \bigr)_{\C^{d}} \log\frac{1}{|z|} dxdy \right)^{1/2}\times \\
\times \left( \iint_{\D}  \bigl( W(z) \overline\partial g(z), 
\overline\partial g(z) \bigr)_{\C^{d}} \log\frac{1}{|z|} dxdy \right)^{1/2}
\end{multline*}
The last term is equivalent to the norm $\|g\|\ci{L^{2}(W)}$ (see 
Lemma \ref{t-norm-equiv}), so by Lemma \ref{l3.1}
$$
|\mathcal I_{2} | \le \|f\|\ci{L^{2}(W)} \cdot \|g\|\ci{L^{2}(W)} 
\cdot
\left\| |z|^{2} \cdot  \| W(z)^{-1/2}\overline\partial W(z) 
W(z)^{-1/2}\| \log \frac1{|z|}  dxdy \right\|_{C}^{1/2} 
$$
Since the measure $\| W(z)^{-1/2}\overline\partial W(z) 
W(z)^{-1/2}\| \log \frac1{|z|}  dxdy$ is a vanishing Carleson measure, 
the Carleson norm $\left\| |z|^{2} \cdot  \| W(z)^{-1/2}\overline\partial W(z) 
W(z)^{-1/2}\| \log \frac1{|z|}  dxdy \right\|_{C}^{1/2} \to 0$ as 
$n\to \infty$. So $|\mathcal I_{2} |\to 0$ as $n\to \infty$. 

To estimate $\mathcal I_{1}$ we pick $r<1$ close to $1$ and split the 
integral into two: $\mathcal I_{1} = \iint_{r\D}\ldots + 
\iint_{\D\setminus r\D} \ldots$. Acting as with $\mathcal I_{2}$ we 
can estimate 
\begin{multline*}
\left| \iint_{X} \ldots \right| \le \\
\le \left(
\iint_{X} \cdot \| W(z)^{-1/2}\overline\partial W(z) 
W(z)^{-1/2}\|^{2}  \cdot 
\bigl( W(z) g(z), 
g(z) \bigr)_{\C^{d}} \log\frac{1}{|z|} dxdy \right)^{1/2}\times \\
\times 
\left( \iint_{X}  \bigl( W(z) \partial\bigl(z^{n} 
f(z)\bigr) , 
\partial\bigl(z^{n} 
f(z)\bigr) \bigr)_{\C^{d}} \log\frac{1}{|z|} dxdy \right)^{1/2} \,,	
\end{multline*}
where $X$ is either $r\D$ or $\D\setminus r\D$. Note that both terms 
are uniformly bounded. 

We can say even more. If $X=r\D$ the second 
term can be made as small as we wish by picking sufficiently large 
$n$. 

Let now $X= \D\setminus r\D$. The measure $\| W(z)^{-1/2}\overline\partial W(z) 
W(z)^{-1/2}\| \log \frac1{|z|}  dxdy$ is a vanishing Carleson 
measure, so for $r$ sufficiently close to $1$ its restriction onto $\D\setminus r\D$ 
has the Carleson norm as small as we want. So by Lemma \ref{l3.1} 
the first term is as small as we want  if $r$ is 
sufficiently close to $1$.

\section{A counterexample to Peller's conjecture.}
\label{c-ex}

In this section we are going to construct a weight $W$, such that 
$W^{-1}\in L^{1}$, $\log W \in\text{VMO}$, but the corresponding 
stationary process is not completely regular (i.e., the  weight $W$ 
does not satisfy any of the conditions 1--6 of Theorem \ref{t-main}). 

Let 
$$
W = U^{*} 
\left(
\begin{array}{cc}
	1 & 0  \\
	0 & \delta(z)
\end{array}
\right)
U \ ,   
\qquad U= \left( 
\begin{array}{cc}
	\cos\alpha & -\sin\alpha  \\
	\sin\alpha  & \cos\alpha
\end{array}\right)\ .
$$
Here 
$$
\delta(e^{it}) = 1/ \log (1/|t|), \qquad -1/4\le t \le 1/4, 
$$
and $\delta$ is a continuous function bounded away from $0$ and $\infty$ 
on the rest of the circle, and 
$$
\alpha(e^{it}) = (t/|t|) \delta(e^{it})^{1/4}, \qquad -1/4\le t \le 1/4, 
$$
and again $\alpha$ is  continuous  on the rest of the circle. 

Then 
$$
\log W = U^{*}
\left(\begin{array}{cc}
	0 & 0  \\
	0 & \log \delta
\end{array}\right) U = 
\left( 
\begin{array}{ll}
	\sin^{2}\alpha \log\delta & \sin\alpha \cos\alpha \log\delta  \\
	 \sin\alpha \cos\alpha \log\delta\quad  & \cos^{2}\alpha \log\delta
\end{array}\right),  
$$
and this matrix clearly belongs to VMO: $\log\delta = \log\log1/|t|$ 
(considered only in a neighborhood of $0$) 
is a ``typical'' unbounded function in VMO, so $\cos^{2}\alpha 
\log\delta\in\text{VMO}$, and all other entries of the matrix are 
continuous. 

Let us now show that the weight $W$ does not even satisfies the 
Muckenhoupt condition $(A_{2})$. Direct computations show that
$$
W= \left( 
\begin{array}{cc}
	\cos^{2}\alpha & -\sin\alpha \cos\alpha  \\
	- \sin\alpha \cos\alpha \quad  & \sin^{2}\alpha 
\end{array}\right) + \delta 
\left(\begin{array}{cc}
	\sin^{2}\alpha & \sin\alpha \cos\alpha  \\
	 \sin\alpha \cos\alpha \quad  & \cos^{2}\alpha 
\end{array}\right)
$$
and 
$$
W^{-1}= \left( 
\begin{array}{cc}
	\cos^{2}\alpha & -\sin\alpha \cos\alpha  \\
	- \sin\alpha \cos\alpha \quad  & \sin^{2}\alpha 
\end{array}\right) + \delta^{-1} 
\left(\begin{array}{cc}
	\sin^{2}\alpha & \sin\alpha \cos\alpha  \\
	 \sin\alpha \cos\alpha \quad  & \cos^{2}\alpha 
\end{array}\right)
$$
If we pick $I$ to be a symmetric arc $[e^{-i\e}, e^{i\e}]$ ($\e>0$ is 
small), then 
off-diagonal entries of $W\ci I$ and $(W^{-1})\ci I$ equal $0$, and 
so we can estimate
$$
W\ci I \ge C
\left(
\begin{array}{cc}
	\cos^{2}\alpha(\e) & 0  \\
	0 & \sin^{2}\alpha(\e)
\end{array}\right) ,
$$
$$
(W^{-1})\ci I \ge 
C\left(
\begin{array}{cc}
	\delta(\e)^{-1}\sin^{2}\alpha(\e) & 0  \\
	0 & \delta(\e)^{-1} \cos^{2}\alpha(\e)
\end{array}\right)\,.
$$
Therefore
$$
\bigl\| [ W\ci I]^{1/2} [(W^{-1}\ci I]^{1/2} \bigr\|\ge C \delta(\e)^{-1} 
\sin\alpha(\e) \cos\alpha(\e) \to \infty\qquad \text{as} \quad \e\to0\,. 
$$


\begin{thebibliography}{30}
\bibitem{Gar} {\sc J.~B.~Garnett,} ``Bounded analytic functions,'' Acad. Press,
NY, 1981. 

\bibitem{HSar} {\sc H. Helson and D. Sarason} Past and Future, {\em 
Math. Scand.,} { \bf 21} (1967), 5--16.

\bibitem{Ib} {\sc I. A. Ibragimov}, Completely regular multidimensional 
stationary processes with discrete time, {\em Proc. Steklov Inst. 
Math.}, {\bf 111} (1970), 269--301.

\bibitem{MW} {\sc P. Masani, N. Wiener}, On bivariate stationary processes and the factorization of matrix -valued functions. {\em Theor. Probability Appl.}, {\bf 4}, (1959), 300-308.  

\bibitem{Pe}{\sc V.V. Peller}, 
 Hankel operators and multivariate stationary processes, \emph{Operator 
 theory: operator algebras and applications, Part 1} (Durham, NH,  
 1988), 357--371, Proc.  Sympos.  Pure Math., 51, Part 1, Amer.  Math.  
 Soc., Providence, RI, 1990.

\bibitem{PKh}{\sc V.V. Peller, S.V. Khruschev} Hankel operators, best approximation, and stationary Gaussian processes, \emph{Russian Math. Surveys} {\bf 37} (1982), 53-124. 

\bibitem{Ro}  
{\sc M. Rosenblum and J. Rovnyak,}
``Hardy classes and operator theory,''  
(Oxford Mathematical Monographs) 
Oxford Science Publications. 
The Clarendon Press, Oxford University Press, New York, 1985.

\bibitem{Roz} {\sc Yu.  A.  Rozanov,} ``Stationary stochastic 
processes,'' Holden-Day, SF, 1967.



\bibitem{Sar} {\sc D. Sarason} An addendum to ``Past and Future'', {\em 
Math. Scand.,} { \bf 30} (1972), 62--64. 

\bibitem{TV1} {\sc S. Treil and A. Volberg}  Wavelets and the angle between 
past and future, \emph{Journal of functional analysis}, {\bf 143}, No. 2, (1997), 269-308.  


\bibitem{TV2} {\sc S. Treil and A. Volberg} Continuous frame decomposition and a 
vector Hunt -- Muckenhoupt -- Wheeden Theorem, 
\emph{ Arkiv f\"{o}r Matematik}, {\bf 35}, No.2, (1997), 363-386. 

\bibitem{TVZ} {\sc S. Treil, A. Volberg and D. Zheng},  Hilbert transform, Toeplitz operators and 
Hankel operators, and invariant $A_{\infty}$ weights, \emph{to appear 
in Rev. Mat. Iberoamericana} 

\bibitem{WM1} {\sc N. Wiener, P. Masani} The prediction theory of multivariate stochastic processes. I. The regularity conditions.
{\em Acta Math.}, {\bf 98}, (1957), 111-150.

\bibitem{WM2}  {\sc N. Wiener, P. Masani}  The prediction theory of multivariate  stochastic processes. II. The linear predicator.{\em Acta Math.}, {\bf 99},
(1958), 93-137. 



\end{thebibliography}
\end{document}